\documentclass[12pt]{amsart} 
\usepackage{amssymb, latexsym, amsthm,epsf}
 
\newtheorem{thm}{Theorem}[section] 

\newtheorem{cor}[thm]{Corollary} 
 
\newtheorem{lemma}[thm]{Lemma}

\newtheorem{ques}[thm]{Question} 
\newtheorem{prop}[thm]{Proposition} 
\newtheorem{rem}[thm]{Remark} 
 
\def\Si{\Sigma} 
 
\def\De{\Delta} 

\def\C{{\mathbb C}} 
 
\def\cL{{\mathcal L}} 
\def\cR{{\mathcal R}} 
\def\F{{\mathbb F}}

\def\P{{\mathbb P}}

\def\Z{{\mathbb Z}}

\def\({\left(} 
\def\){\right)} 

\long\def\forget#1\forgotten{}

\newif \iffurther 
\furtherfalse 
 
\newif \iffurther 
\furtherfalse 
 
\newif\ifXY 
\XYfalse    
 
\ifXY 
\usepackage{xy} 
\fi 
\ifXY 
\xyoption{all} 
\fi 
 
\begin{document}

\title[A connection between fundamental groups]{On the connection between
  affine and projective fundamental groups of line arrangements and curves} 
 
\author[David Garber]{David Garber$^1$} 
 
\stepcounter{footnote} 
\footnotetext{Partially supported by the Lady Davis postdoctoral fellowship.} 
 
\address{Einstein Institute of Mathematics, The Hebrew University, Givat Ram, 91904 Jerusalem, Israel} 
\email{garber@math.huji.ac.il} 
  
\date{\today} 
 
\begin{abstract} 
In this note we prove a decomposition related to the affine
fundamental group and the projective fundamental group of 
a line arrangement and a reducible curve with a line component.
We give some applications to this result. 
\end{abstract} 
 
\maketitle 
 
\section{Introduction} 

The fundamental group of complements of plane curves is a very
important topological invariant with many different applications.
There are two such invariants: the {\it affine fundamental group}
of a plane curve, which is the fundamental group of its affine complement, 
and its {\it projective fundamental group}, 
which is the fundamental group of its projective complement. 

Oka \cite{O} has proved the following interesting result, which shed new light
on the connection between these two fundamental groups:

\begin{thm}[Oka] \label{oka}
Let $C$ be a curve in $\C\P ^2$ and let $L$ be 
a general line to $C$, {\it i.e.} $L$ intersects $C$ in only simple points.
Then, we have a central extension:
$$1 \to \Z \to \pi _1(\C\P ^2 -(C \cup L)) \to \pi _1(\C\P^2-C) \to 1$$ 
\end{thm}

Since $L$ is a general line to $C$, then:
$$\pi _1(\C\P ^2 -(C \cup L)) =\pi _1(\C^2 -C).$$
Hence, we get the following interesting connection 
between the two fundamental groups:
$$1 \to \Z \to \pi _1(\C^2 -C) \to \pi _1(\C\P^2-C) \to 1$$ 
  
\medskip

A natural question is: 
\begin{ques}
Under which conditions does this short exact sequence split? 
Notice that when it does, we have the following decomposition:
$$\pi _1(\C^2 -C) \cong \pi _1(\C\P^2-C) \oplus \Z$$ 
\end{ques}
 
A {\it real line arrangement} in $\C^2$ is a finite union of copies of $\C$ in $\C^2$, whose
equations can be written by real coefficients. 
Some families of real line arrangements were already proved to satisfy this
condition, see \cite{GaTe}.

Here, we show that:
\begin{thm}\label{main}
If $\cL$ is a real line arrangement, then such a decomposition holds:
$$\pi _1(\C^2 -\cL) \cong \pi _1(\C\P^2-\cL) \oplus \Z$$ 
\end{thm}

Actually, one can see that the same argument holds for arbitrary 
line arrangements (see Theorem \ref{complex_la}). 
Moreover, we give a different condition for this decomposition to hold: 
If $C$ is a plane curve with a line component, i.e. $C = C' \cup L$
where $L$ is a line, we have that
$\pi_1(\C^2 -C) =\pi_1(\C\P^2 -C)\oplus \Z$ too.  

These results can be applied to the computation of the affine fundamental groups
of line arrangements, since the projective fundamental group is an
easier object to deal with than the affine fundamental group.

\medskip

The paper is organized as follows.
In Section \ref{proof} we prove Theorem \ref{main}.
The other condition is discussed in Section \ref{fg},
and in the last section we give some applications.

\section{Proof of Theorem \ref{main}}\label{proof}

In this section, we prove Theorem \ref{main}. 

\begin{proof}[Proof of Theorem \ref{main}]
Let $\cL$ be a real line arrangement with $n$ lines. 
Let $L$ be an arbitrary line which intersects $\cL$ transversally. By the following remark,
one can reduce the proof to the case where $L$ is a line which intersects
transversally all the lines in $\cL$  and all the intersection points of $L$ with
lines in $\cL$ are to the left of all the intersection points of $\cL$ (see Figure
\ref{arrangement} for such an example, where the arrangement $\cL$ 
consists of $L_1,L_2,L_3$ and $L_4$).

\begin{figure}[h]
\epsfysize=5cm 
\centerline{\epsfbox{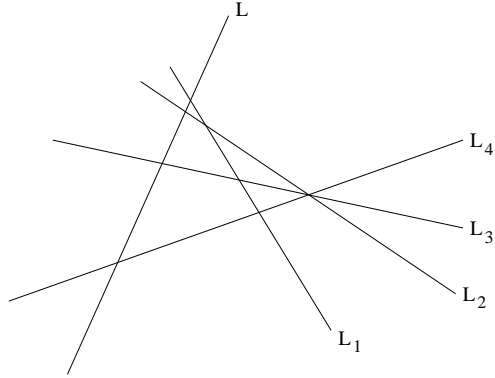}} 
\caption{An example}\label{arrangement} 
\end{figure}                                         
    
\begin{rem}\label{line_infty}
We have proved in \cite{GTV1} that in a real line arrangement, if a line crosses 
a multiple intersection point from one side to its other side (see Figure \ref{pass_mul}), 
the fundamental groups remain unchanged  
(see \cite[Theorem 4.13]{GTV1} for this property of the action 
$\stackrel{\Delta}{=}$). 

\begin{figure}[h]
\epsfysize=2cm 
\centerline{\epsfbox{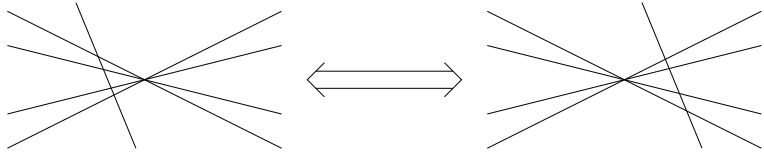}} 
\caption{A line crosses a multiple intersection point.}\label{pass_mul} 
\end{figure}                                         

By this argument, we can start with any transversal line $L$ as  
``the line at infinity''. Then, by using this property repeatedly, 
we can push this line over all 
the intersection points of the arrangement $\cL$, without changing the corresponding 
fundamental group. This process will be terminated when all the intersection points of $L$ with 
$\cL$ are placed to the left of all the intersection points of $\cL$, and this is the reduced 
case. 
\end{rem}

By this remark, we continue the proof of Theorem \ref{main} for the reduced case, 
where all the intersection points of $L$ with 
$\cL$ are placed to the left of all the intersection points of $\cL$.
We compute presentations for $\pi_1(\C\P^2 -\cL)$ and 
$\pi_1(\C\P^2 -(\cL \cup L))$
by braid monodromy techniques (the Moishezon-Teicher algorithm) and 
the van Kampen Theorem.
The original techniques are introduced in \cite{BGT1} 
and \cite{vK} respectively. Shorter presentations of these 
techniques can be found in \cite{GaTe} and \cite{GTV1}. 

We first have to compute the lists of Lefschetz pairs of the arrangements, 
which are the pairs of indices of the components intersected 
at the intersection 
point, where we numerate the components locally near the point 
(see \cite{GaTe}).
Since $\cL$ is an arbitrary real line arrangement, its list of Lefschetz pairs is
$$([a_1,b_1], \cdots ,[a_k,b_k]),$$ 
where $k$ is the number of intersection points in 
$\cL$. If we assume that the additional line $L$ crosses all the lines of $\cL$ 
transversally to the left of all the intersection points in $\cL$, the list 
of Lefschetz pairs of $\cL \cup L$ is 
$$([a_1,b_1], \cdots ,[a_k,b_k],[n,n+1], \cdots , [1,2])$$  

By the braid monodromy techniques and the van Kampen Theorem, 
the group $\pi_1(\C\P^2 -\cL)$ has $n$ generators $\{ x_1, \cdots, x_n \}$ and
$\pi_1(\C\P^2 -(\cL \cup L))$ has $n+1$ generators $\{ x_1, \cdots, x_n,x_{n+1} \}$.
Moreover, the first $k$ relations of 
$\pi_1(\C\P^2 -(\cL \cup L))$ are equal to the $k$ relations of $\pi_1(\C\P^2 -\cL)$.
Let us denote this set of relations by $\cR$. 

Now, we have to find out the relations induced by the $n$ intersection 
points of the line $L$ with the arrangement $\cL$. Moreover, we have 
to add at last the appropriate projective relations.

First, we compute the relations induced by the $n$ intersection 
points of the line $L$ with the arrangement $\cL$. The Lefschetz pair 
of the $(k+1)$th point (which is the first intersection point of $L$ with $\cL$) 
is $[n,n+1]$. Hence, its corresponding initial skeleton is shown in Figure 
\ref{init_skel_k_1}.

\begin{figure}[h]
\epsfxsize=7cm 
\centerline{\epsfbox{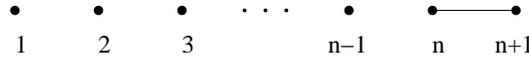}} 
\caption{Initial skeleton of the $(k+1)$th point}\label{init_skel_k_1} 
\end{figure}                                         
 
Since the list of pairs  $[a_1,b_1], \cdots ,[a_k,b_k]$ is induced by a 
line arrangement, it is easy to see that:
$$ \De \langle a_k,b_k \rangle \cdots \De  \langle a_1,b_1 \rangle = \De \langle 1,n \rangle,$$
where $\De \langle 1,n \rangle$ is the generalized half-twist on the segment $[1,n]$.

Therefore, applying this braid on the initial skeleton yields the resulting skeleton 
for this point which is presented in Figure \ref{final_skel_k_1}.

\begin{figure}[h]
\epsfxsize=7cm 
\centerline{\epsfbox{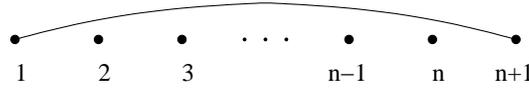}} 
\caption{Final skeleton of the $(k+1)$th point}\label{final_skel_k_1} 
\end{figure}                                         

By the van Kampen Theorem, the corresponding relation is:
$$[x_n x_{n-1} \cdots x_2 x_1 x_2^{-1} \cdots x_n^{-1},x_{n+1}] = 1$$     

Let $P_i$ be the $i$th intersection point, where $k+2 \leq i \leq k+n$.
Its corresponding Lefschetz pair is $[n-(i-k)+1,n-(i-k)+2]$, and hence
its initial skeleton is shown in Figure \ref{init_skel_k_i}.

\begin{figure}[h]
\epsfxsize=8cm 
\centerline{\epsfbox{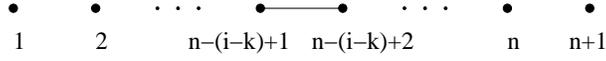}} 
\caption{Initial skeleton of the $i$th point}\label{init_skel_k_i} 
\end{figure}

Now, we first have to apply on it the following braid:
$$\De \langle  n-(i-k)+2,n-(i-k)+3 \rangle \cdots \De \langle n,n+1 \rangle$$
and afterwards we have to apply on the resulting skeleton $\De \langle 1,n \rangle$, 
which equals to 
$$\De \langle a_k,b_k \rangle \cdots \De  \langle a_1,b_1 \rangle$$
as before. Hence, Figure \ref{final_skel_k_i} presents the resulting skeleton 
for the $i$th point.

\begin{figure}[h]
\epsfxsize=9cm 
\centerline{\epsfbox{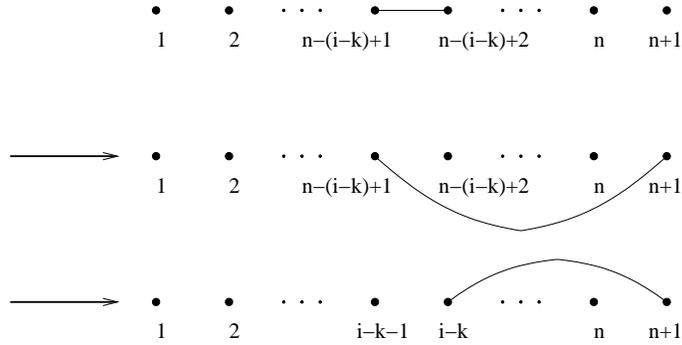}} 
\caption{Computing the final skeleton of the $i$th point}\label{final_skel_k_i} 
\end{figure}                                         

Again, by the van Kampen Theorem, the corresponding relation is:
$$[x_n x_{n-1} \cdots x_{i-k+1} x_{i-k} x_{i-k+1}^{-1} \cdots x_n^{-1},x_{n+1}] = 1$$     
To summarize, we get that the set of relations induced by the intersection points
of the additional line $L$ is:
$$\{ [x_n x_{n-1} \cdots x_{i-k+1} x_{i-k} x_{i-k+1}^{-1} \cdots x_n^{-1},x_{n+1}] = 1 \quad | \quad 1 \leq i-k \leq n \}$$
One can easily see, by a sequence of substitutions, that actually 
this set of relations is equal to the following set:
$$\{ [x_i ,x_{n+1}] = 1 \quad | \quad 1 \leq i \leq n \}$$

Hence, we have the following presentations:

$$\pi_1(\C\P^2-\cL) = \langle x_1,\cdots, x_n \quad | \quad \cR \ ;\ x_n x_{n-1} \cdots x_1=1  \rangle$$
$$\pi_1(\C\P^2-(\cL \cup L)) = 
\left\langle 
\begin{array}{c|c}
x_1,\cdots, x_n,x_{n+1} & \cR \ ;\ [x_i,x_{n+1}]=1\ ,\ 1 \leq i \leq n\ ; \\
 & x_{n+1} x_n \cdots x_1=1
\end{array}
\right\rangle$$

It remains to show that these presentations imply that: 
$$\pi_1(\C\P^2-(\cL \cup L)) = \Z \oplus \pi_1(\C\P^2-\cL)$$   

Denote by $\cR(x_1 \Leftarrow w)$ the set of relations $\cR$ after we substitute anywhere
the generator $x_1$ by an element $w$ in the corresponding group.

Hence, using the projective relations in both fundamental groups, we have the 
following new presentations:
$$\pi_1(\C\P^2-\cL) = \langle x_2,\cdots, x_n \quad | \quad \cR(x_1 \Leftarrow x_2^{-1} \cdots x_n^{-1}) \rangle$$
$$\pi_1(\C\P^2-(\cL \cup L)) = \left\langle 
\begin{array}{c|c}
x_2,\cdots, x_n,x_{n+1} &  \cR(x_1 \Leftarrow x_2^{-1} \cdots x_n^{-1} x_{n+1}^{-1});\\
 & [x_i,x_{n+1}]=1, 2 \leq i \leq n  
\end{array}
\right\rangle$$

Since all the relations in $\cR$ are either commutative relations 
or cyclic relations of the type 
$$w_m w_{m-1} \cdots w_1 = w_{m-1} \cdots w_1 w_m = \cdots = w_1 w_m \cdots w_2,$$
where $w_i$ are words in the generators $x_1, \cdots , x_n$, then the sum of powers
of the generator $x_1$ or its inverse are equal in any part 
of the relation. Hence, using the fact that $x_{n+1}$ commutes with all the 
other generators, we actually can cancel $x_{n+1}$ from the set of relations
$$\cR(x_1 \Leftarrow x_2^{-1} \cdots x_n^{-1} x_{n+1}^{-1}),$$
and therefore we have:
$$\cR(x_1 \Leftarrow x_2^{-1} \cdots x_n^{-1} x_{n+1}^{-1})=\cR(x_1 \Leftarrow x_2^{-1} \cdots x_n^{-1}).$$

By this argument, we now have: 
\begin{eqnarray*}
\pi_1(\C\P^2-(\cL \cup L)) & \cong & 
\left\langle 
\begin{array}{c|c}
x_2,\cdots, x_n,x_{n+1} & \cR(x_1 \Leftarrow x_2^{-1} \cdots x_n^{-1} x_{n+1}^{-1}) ;\\
 & [x_i,x_{n+1}]=1, 2 \leq i \leq n 
\end{array}
\right\rangle \\
 & \cong & 
\left\langle 
\begin{array}{c|c}
x_2,\cdots, x_n,x_{n+1} &  \cR(x_1 \Leftarrow x_2^{-1} \cdots x_n^{-1});\\
 & [x_i,x_{n+1}]=1, 2 \leq i \leq n  
\end{array}
\right\rangle \\
 & \cong & \langle x_2,\cdots, x_n \quad | \quad \cR(x_1 \Leftarrow x_2^{-1} \cdots x_n^{-1})\rangle \oplus \langle x_{n+1} \rangle \\
 & \cong & \pi_1(\C\P^2-\cL) \oplus \Z  \\
\end{eqnarray*}

Since $L$ is transversal to $\cL$, this is equivalent to:
$$\pi_1(\C^2-\cL) \cong \pi_1(\C\P^2-\cL) \oplus \Z$$
and we are done.
\end{proof}

\begin{rem}
By the projective relation in $\pi_1(\C\P^2-(\cL \cup L))$, we have that:
$$x_{n+1} = (x_n \cdots x_1)^{-1}.$$
Hence, we have that the element which corresponds to the loop around all the lines
is a central element in $\pi_1(\C^2-\cL)=\pi_1(\C\P^2-(\cL \cup L))$, as was shown 
by Oka \cite{O} too. 
\end{rem}

Now, one can see that a similar argument holds for complex 
line arrangements too. Hence, one has the following generalization:

\begin{thm}\label{complex_la}
Let $\cL$ be an arbitrary line arrangement. Then the following 
decomposition holds:
$$\pi _1(\C^2 -\cL) \cong \pi _1(\C\P^2-\cL) \oplus \Z$$ 
\end{thm}

\section{A different condition for a decomposition}\label{fg} 

In this section we will show that the decomposition holds for 
any reducible curve with a line component. 

\medskip

Oka and Sakamoto \cite{OS} have proved the following decomposition  
theorem concerning the affine fundamental group: 
\begin{thm}[Oka-Sakamoto]\label{OS} 
Let $C_1$ and $C_2$ be algebraic plane curves in $\C ^2$.  
Assume that the intersection $C_1 \cap C_2$  
consists of distinct $d_1 \cdot  d_2$ points, where $d_i \ (i=1,2)$ are the  
respective degrees of $C_1$ and $C_2$. Then: 
$$\pi _1 (\C ^2 - (C_1 \cup C_2)) \cong \pi _1 (\C ^2 -C_1) \oplus \pi _1 (\C ^2 -C_2)$$ 
\end{thm} 
 
Here is a parallel version for the projective fundamental group: 
 
\begin{lemma} 
Let $C_1$ and $C_2$ be plane curves in $\C\P ^2$ with respective degrees $d_1$  
and $d_2$. Let $C_i'$ be the affine part of $C_i$ with respect to a
projective line $L$. 
If $C_1' \cap C_2'$ consists of $d_1 d_2$ distinct simple points, then:    
$$\pi _1 (\C\P ^2 - (C_1 \cup C_2 \cup L)) \cong \pi _1 (\C\P ^2 - (C_1 \cup L)) \oplus \pi _1 (\C\P ^2 -(C_2 \cup L)).$$  
\end{lemma} 
 
\begin{proof} 
When $L$ becomes the line at infinity,  
the curve $C_1 \cup C_2$ becomes $C_1' \cup C_2'$. Also, it is easy to see  
that if $C'$ denotes the affine part of $C$ with respect to $L$, then 
$$\pi _1(\C\P ^2 -(C \cup L)) \cong \pi _1(\C^2-C').$$ 
 
Now, we can compute: 
\begin{eqnarray} 
\pi _1 (\C\P ^2 - (C_1 \cup C_2 \cup L)) & \cong &  \pi_1 (\C ^2 - (C_1 ' \cup  
C_2')) \nonumber \\ 
   & \stackrel{\rm{[7]}}{\cong} & \pi _1 (\C ^2-C_1') \oplus \pi _1(\C ^2- 
 C_2 ') \nonumber \\ 
   & \cong & \pi _1(\C\P ^2- (C_1 \cup L)) \oplus \pi _1(\C\P ^2-(C_2 \cup L)) 
 \nonumber 
\end{eqnarray} 
\end{proof} 
  
Hence, we get:

\begin{cor} 
Let $C$ be a plane curve of degree $d$ in $\C\P ^2$ and 
let $\mathcal L$ be a line arrangement 
consisting of $k$ lines meeting in a point outside $C$. 
Let $L$ be an additional line which intersects both  
$C$ and $\mathcal L$ transversally. 
If $C \cap \mathcal L$ consists of $d k$ distinct simple points, then:    
$$\pi _1 (\C\P ^2 - (C \cup \mathcal L \cup L)) \cong \pi _1 (\C\P ^2 - (C \cup L)) \oplus \Z \oplus \F_{k-1},$$ 
where $\F_k$ is the free group with $k$ generators.
In particular, if $k=1$, i.e. $\mathcal L$ consists of one line $L_1$, then: 
$$\pi _1 (\C\P ^2 - (C \cup L_1 \cup L)) \cong \pi _1 (\C\P ^2 - (C \cup L)) \oplus \Z$$ 
\end{cor} 
 
\begin{proof} 
Substituting $C$ instead of $C_1$ and $\mathcal L$ instead of $C_2$ in the last lemma, 
yields the following equation: 
$$\pi _1 (\C\P ^2 - (C \cup \mathcal L \cup L)) \cong \pi _1 (\C\P ^2 - (C \cup L)) \oplus \pi _1 (\C\P ^2 -(\mathcal L \cup L)).$$ 
Since $\pi _1 (\C\P ^2 -(\mathcal L \cup L)) = \Z \oplus \F_{k-1}$ (see for example \cite{GaTe}), 
the Corollary is proved. Fixing $k=1$ yields the particular case. 
\end{proof} 

Hence we got the following result:
\begin{prop}\label{line_component}
Let $C$ be a plane curve with a line component, i.e. $C = C' \cup L$
where $L$ is a line that intersects $C'$ transversally. Then:
$$\pi_1(\C^2 -C) =\pi_1(\C\P^2 -C)\oplus \Z$$  
\end{prop}

\begin{proof} 
Using the case $k=1$ in the last corollary, we get that:
$$\pi _1 (\C\P ^2 - (C' \cup L_1 \cup L)) \cong \pi _1 (\C\P ^2 - (C' \cup L)) \oplus \Z$$ 
But since $L_1$ intersects $C' \cup L$ transversally, we have that
$$\pi _1 (\C\P ^2 - (C' \cup L_1 \cup L)) \cong \pi _1 (\C ^2 - (C' \cup
L))$$
So we get:
$$\pi _1 (\C ^2 - (C' \cup L)) \cong \pi _1 (\C\P ^2 - (C' \cup L)) \oplus \Z$$ 
\end{proof} 

\begin{rem}
By a similar argument to the argument we have used to prove Theorem \ref{main},
and by the fact that the line at infinity is transversal to $C \cup L$, 
we can generalize Proposition \ref{line_component} to the case that 
$L$ does not have to be transversal to $C$. 
\end{rem}

\section{Some applications}\label{applications}

In this section, we present some immediate applications of the main result.

\subsection{Structure of the affine fundamental groups for some real arrangements}

Fan \cite{Fa} proved the following result:
\begin{prop}[Fan]
Let $\Si$ be an arrangement of $n$ lines and $S=\{a_1, \cdots, a_k\} $
be the set of all singularities of $\Si$ with multiplicity $\geq 3$. 
Suppose that the subgraph of $\Si$, which contains only the higher singularities 
(i.e. with multiplicity $\geq 3 $) and their edges, is a union of trees. Then:
$$\pi _1 (\C \P ^2 - \Si) \cong \Z ^r \oplus \F ^{m(a_1)-1} \oplus \cdots \oplus \F ^{m(a_k)-1}$$
where $r=n+k-1-m(a_1)- \cdots - m(a_k)$ and $m(a_i)$ is the multiplicity
of the intersection point $a_i$.
\end{prop}

Using Theorem \ref{main}, we have the following easy consequence:

\begin{cor}
Under the same conditions of the last proposition, we have that:
$$\pi _1 (\C ^2 - \Si) \cong \Z ^{r+1} \oplus \F ^{m(a_1)-1} \oplus \cdots \oplus \F ^{m(a_k)-1}$$
where $r=n+k-1-m(a_1)- \cdots - m(a_k)$.
\end{cor}

This generalizes Theorem 5.3 of \cite{GaTe}.

\subsection{The connection between the invariants of real arrangements}
In \cite{GTV1} and \cite{GTV2} we address the following question 
about the invariants of real arrangements:

\begin{ques}
Does the incidence lattice determine the affine fundamental group and 
the projective fundamental group of real line arrangements?
\end{ques}

For complex line arrangements, the answer for this question is negative 
due to an example of Rybnikov (see \cite{Ry}).

\medskip

Our main result is (see \cite{GTV2}):
\begin{thm}
For real line arrangements with up to $8$ lines, the incidence lattice 
does determine the affine fundamental group and the projective fundamental group.
\end{thm}

By Theorem \ref{main}, we actually reduce this question to the projective 
case only, since we have:

\begin{cor}
Let $\cL_1$ and $\cL_2$ be two real line arrangements such that:
$$ \pi_1(\C\P^2- \cL_1) \cong \pi_1(\C\P^2- \cL_2)$$
Then:
$$ \pi_1(\C^2- \cL_1) \cong \pi_1(\C^2- \cL_2)$$
\end{cor}
 
This simplifies much the computational aspect of this question, since during
the computations we have performed in \cite{GTV2}, it was turned out that 
the comparison of two projective fundamental groups was much easier 
than the comparison of the corresponding affine fundamental groups. 
The main reason for this difference is probably that the affine 
fundamental group has one more generator than the projective fundamental 
group (this generator is cancelled by the projective relation which implies that 
the multiplication of all the generators is the unit element). 
 
\section*{Acknowledgment} 
I wish to thank Mutsuo Oka and Pho Duc Tai for fruitful discussions.
I also thank Ehud De-Shalit and Einstein Institute of Mathematics 
for hosting my stay.

\end{document}